\documentclass[12pt,twoside,a4paper]{article}
\usepackage[english]{babel}
\usepackage{amssymb}
\usepackage[T1]{fontenc}
\usepackage{amsmath,epsf}
\newtheorem{lemma}{Lemma}

\newtheorem{Prop}{Proposition}

\newcommand{\A}{\mathcal{A}}
\newcommand{\C}{\mathcal{C}}
\newcommand{\F}{\mathcal{F}}
\newcommand{\R}{\mathbb{R}}
\newcommand{\Z}{\mathbb{Z}}
\newcommand{\N}{\mathbb{N}}

\newcommand{\B}{\mathcal{B}}
\newcommand{\I}{\mathcal{I}}
\newcommand{\U}{\mathcal{U}}
\newcommand{\V}{\mathcal{V}}

\renewcommand{\P}{\mathbb{P}}
\newcommand{\ds}{\displaystyle}
\newcounter{tictac}
\newenvironment{fleuveA}{
   \begin{list}{$\textbf{A\arabic{tictac}}$) }{\usecounter{tictac}
\leftmargin 1cm\labelwidth 2em}}{\end{list}}
\def\1{\,\rlap{\mbox{\small\rm 1}}\kern.15em 1}
\def\ind#1{\1_{#1}}
\def\build#1_#2^#3{\mathrel{\mathop{\kern 0pt#1}\limits_{#2}^{#3}}}
\def\tend#1#2{\build\hbox to 12mm{\rightarrowfill}_{#1\rightarrow #2}^{a.s.}}

\def\converge#1#2#3{\build\hbox to
15mm{\rightarrowfill}_{#1\rightarrow #2}^{\hbox{\scriptsize #3}}}
\begin{document}
\title{Asymptotic normality of kernel estimates in a regression model for random fields}
\author{Mohamed EL MACHKOURI, Radu STOICA}
\maketitle

{\renewcommand\abstractname{Abstract}
\begin{abstract}
We establish the asymptotic normality of the regression estimator in a
fixed-design setting when the errors are given by a field of dependent
random variables. The result applies to martingale-difference or
strongly mixing random fields. On this basis, a statistical test that can
be applied to image analysis is also presented.
{\em Short title:} Asymptotic normality of kernel
estimates.
\end{abstract}
%-------------------------------------------------- Introduction ----------------------------------------------------------
\thispagestyle{empty}
\section{Introduction and notations}
Our aim in this paper is to establish the asymptotic normality of
a regression estimator in a fixed-design setting when the errors
are given by a stationary field of random variables which show
spatial interaction. Let $\Z^{d}$, $d\geq 1$ denote the integer
lattice points in the $d$-dimensional Euclidean space. By a
stationary random field we mean any family
$(\varepsilon_{k})_{k\in\Z^{d}}$ of real-valued random variables
defined on a probability space $(\Omega, \F, \P)$ such that for
any $(k,n)\in\Z^{d}\times\N^{\ast}$ and any $(i_{1},...,i_{n})\in
(\Z^{d})^{n}$, the random vectors
$(\varepsilon_{i_{1}},...,\varepsilon_{i_{n}})$ and
$(\varepsilon_{i_{1}+k},...,\varepsilon_{i_{n}+k})$ have the same
law. The regression model which we are interested in is
\begin{equation}\label{model}
Y_{i}=g(i/n)+\varepsilon_{i},\quad i\in\Lambda_{n}=\{1,...,n\}^{d}
\end{equation}
where $g$ is an unknown smooth function and
$(\varepsilon_{i})_{i\in\Z^{d}}$ is a zero mean and square-integrable
stationary random field. Let $K$ be a probability kernel defined on $\R^{d}$
and $(h_{n})_{n\geq 1}$ a sequence of positive numbers which
converges to zero and which satisfies $(nh_{n})_{n\geq 1}$ goes to
infinity. We estimate the function $g$ by the kernel-type
estimator $g_{n}$ defined for any $x$ in $[0,1]^{d}$ by
\begin{equation}\label{def-g_n}
g_{n}(x)=\frac{\ds{\sum_{i\in\Lambda_{n}}}\,Y_{i}\,
K\left(\frac{x-i/n}{h_{n}}\right)}{\ds{\sum_{i\in\Lambda_{n}}}\,
K\left(\frac{x-i/n}{h_{n}}\right)}.
\end{equation}
In a previous paper, El Machkouri \cite{Elmachkouri-SISP} obtained
strong convergence of the estimator $g_{n}(x)$ with optimal rate. However, most of
existing theoretical nonparametric results for dependent random
variables pertain to time series (see Bosq \cite{Bosq-livre}) and
relatively few generalisations to the spatial domain are
available. Key references on this topic are Biau \cite{Biau-2003},
Carbon et al. \cite{Carbon-Hallin-Tran}, Carbon et al.
\cite{Carbon-Tran-Wu}, Hallin et al. \cite{Hallin-Lu-Tran-2001},
\cite{Hallin-Lu-Tran-2004a}, Tran \cite{Tran}, Tran and Yakowitz
\cite{Tran-Yakowitz} and Yao \cite{Yao} who have investigated
nonparametric density estimation for random fields and Altman
\cite{Altman}, Biau and Cadre \cite{Biau-Cadre}, Hallin et al.
\cite{Hallin-Lu-Tran-2004b} and Lu and Chen \cite{Lu-Chen-2002},
\cite{Lu-Chen-2004} who have studied spatial prediction and
spatial regression estimation.\\
Let $\mu$ be the law of the stationary real random field
$(\varepsilon_{k})_{k\in\Z^{d}}$ and consider the projection $f$
from $\R^{\Z^{d}}$ to $\R$ defined by $f(\omega)=\omega_{0}$ and
the family of translation operators $(T^{k})_{k\in\Z^{d}}$ from
$\R^{\Z^{d}}$ to $\R^{\Z^{d}}$ defined by
$(T^{k}(\omega))_{i}=\omega_{i+k}$ for any $k\in\Z^{d}$ and any
$\omega$ in $\R^{\Z^{d}}$. Denote by $\B$ the Borel
$\sigma$-algebra of $\R$. The random field $(f\circ
T^{k})_{k\in\Z^{d}}$ defined on the probability space
$(\R^{\Z^{d}}, \B^{\Z^{d}}, \mu)$ is stationary with the same law
as $(\varepsilon_{k})_{k\in\Z^{d}}$, hence, without loss of
generality, one can suppose that $(\Omega, \F, \P)=(\R^{\Z^{d}},
\B^{\Z^{d}}, \mu)$ and $\varepsilon_{k}=f\circ T^{k}$. An element
$A$ of $\F$ is said to be invariant if $T^{k}(A)=A$ for any
$k\in\Z^{d}$. We denote by $\I$ the $\sigma$-algebra of all
measurable invariant sets. On the lattice $\Z^{d}$ we define the
lexicographic order as follows: if $i=(i_{1},...,i_{d})$ and
$j=(j_{1},...,j_{d})$ are distinct elements of $\Z^{d}$, the
notation $i<_{lex}j$ means that either $i_{1}<j_{1}$ or for some
$p$ in $\{2,3,...,d\}$, $i_{p}<j_{p}$ and $i_{q}=j_{q}$ for $1\leq
q<p$. Let the sets
$\{V_{i}^{k}\,;\,i\in\Z^{d}\,,\,k\in\N^{\ast}\}$ be defined as
follows:
$$
V_{i}^{1}=\{j\in\Z^{d}\,;\,j<_{lex}i\},
$$
and for $k\geq 2$
$$
V_{i}^{k}=V_{i}^{1}\cap\{j\in\Z^{d}\,;\,\vert i-j\vert\geq
k\}\quad\textrm{where}\quad \vert i-j\vert=\max_{1\leq l\leq
d}\vert i_{l}-j_{l}\vert.
$$
For any subset $\Gamma$ of $\Z^{d}$ define
$\F_{\Gamma}=\sigma(\varepsilon_{i}\,;\,i\in\Gamma)$ and set
$$
E_{\vert k\vert}(\varepsilon_{i})=E(\varepsilon_{i}\vert\F_{V_{i}^{\vert
k\vert}}),\quad k\in V_{i}^{1}.
$$
Note that Dedecker \cite{JD-tcl} established the central limit
theorem for any stationary square-integrable random field
$(\varepsilon_{k})_{k\in\Z^{d}}$ which satisfies the condition
\begin{equation}\label{proj}
\sum_{k\in V_{0}^{1}}\|\varepsilon_{k}E_{\vert
k\vert}(\varepsilon_{0})\|_{1}<\infty.
\end{equation}
A real random field $(X_{k})_{k\in\Z^{d}}$ is said to be a
martingale-difference random field if for any $m$ in $\Z^{d},\,\,E(\,X_{m}\,\vert\,\sigma(\,X_{k}\,;\,k<_{lex}m\,)\,)=0$ a.s.
The condition ($\ref{proj}$) is satisfied by martingale-difference random fields.
Nahapetian and Petrosian \cite{Nahapetian-Petrosian} defined a large class of Gibbs
random fields $(\xi_{k})_{k\in\Z^{d}}$  satisfying the stronger martingale-difference property:
$E(\,\xi_{m}\,\vert\,\sigma(\,\xi_{k}\,;\,k\neq m\,)\,)=0$ a.s. for any $m$ in $\Z^{d}$.
Moreover, for these models, phase transition may occur (see \cite{Martirosian},\cite{Nahapetian}). \\
\\
Given two sub-$\sigma$-algebras $\U$ and $\V$, different measures of their dependence have been considered
in the literature. We are interested by one of them. The strong mixing (or $\alpha$-mixing) coefficient has been
introduced by Rosenblatt \cite{Ros} and is defined by
$$
\alpha(\U,\V)=\sup\{\vert\P(U\cap V)-\P(U)\P(V)\vert,\,U\in\U
,\,V\in\V\}.
$$
Denote by $\sharp\Gamma$ the cardinality of any subset $\Gamma$ of
$\Z^{d}$. In the sequel, we shall use the following non-uniform
mixing coefficients defined for any $(k,l,n)$ in
$(\N^{\ast}\cup\{\infty\})^{2}\times\N$ by
$$
\alpha_{k,l}(n)=\sup\,\{\alpha(\F_{\Gamma_{1}},\F_{\Gamma_{2}}),\,
\sharp\Gamma_{1}\leq k,\, \sharp\Gamma_{2}\leq l,\,
\rho(\Gamma_{1},\Gamma_{2})\geq n\},
$$
where the distance $\rho$ is defined by
$\rho(\Gamma_{1},\Gamma_{2})=\min\{\vert
i-j\vert,\,i\in\Gamma_{1},\,j\in\Gamma_{2}\}$. We say that the
random field $(\varepsilon_{k})_{k\in\Z^{d}}$ is strongly mixing (or $\alpha$-mixing)
if there exists a pair $(k,l)$ in $(\N^{\ast}\cup\{\infty\})^{2}$
such that $\lim_{n\to \infty}\alpha_{k,l}(n)=0$. \\
The condition ($\ref{proj}$) is satisfied by strongly mixing random fields. For example,
one can construct stationary Gaussian random fields with a sufficiently large polynomial
decay of correlation such that (\ref{condition-sm-strong-mixing}) holds (\cite{Doukhan}, p. 59, Corollary 2).
\section{Main results}
First, we recall the concept of stability introduced by R\'enyi
\cite{Renyi}.\\
\vspace{-0.2cm}
\\
\textbf{Definition.} {\em Let $(X_{n})_{n\geq 0}$ be a sequence
of real random variables and let $X$ be defined on some extension
of the underlying probability space $(\Omega, \A, \P)$. Let $\U$
be a sub-$\sigma$-algebra of $\A$. Then $(X_{n})_{n\geq 0}$ is said to
converge $\U$-stably to $X$ if for any continuous bounded function
$\varphi$ and any bounded and $\U$-measurable variable $Z$ we have
$\lim_{n\to\infty}E\left(\varphi(X_{n})Z\right)=E\left(\varphi(X)Z\right)$.}\\
\vspace{-0.2cm}
\\
For any $B>0$, we denote by $\C^{1}(B)$ the set of
real functions $f$ continuously differentiable on $[0,1]^{d}$ such
that
$$\sup_{x\in [0,1]^{d}}
\max_{\alpha\in\mathcal{M}}\vert D_{\alpha}(f)(x)\vert\leq B,
$$
where
$$
D_{\alpha}(f)=\frac{\partial^{\hat{\alpha}}f}{\partial
x_{1}^{\alpha_{1}}...\,\partial
x_{d}^{\alpha_{d}}}\quad\textrm{and}
\quad\mathcal{M}=\{\alpha=(\alpha_{i})_{i}
\in\N^{d}\,;\,\hat{\alpha}=\sum_{i=1}^{d}\alpha_{i}\leq 1\}.
$$
In the sequel we denote $\|x\|=\max_{1\leq k\leq d}\vert
x_{k}\vert$ for any $x=(x_{1},...,x_{d})\in[0,1]^d$. We make the following assumptions on
the regression function $g$ and the probability kernel K:
\begin{fleuveA}
\item The probability kernel $K$ fulfils $\int K(u)\,du=1$ and $\int K^2(u)\,du<\infty$. $K$ is also symmetric, non-negative,
supported by $[-1,1]^{d}$ and satisfies a Lipschitz condition
$\vert K(x)-K(y)\vert\leq r\|x-y\|$ for any $x,y\in[-1,1]^{d}$
and some $r>0$. In addition there exists $c,C>0$ such that
$c\leq K(x)\leq C$ for any $x\in[-1,1]^{d}$.
\item There exists $B>0$ such that $g$ belongs to $\C^{1}(B)$.
\end{fleuveA}
We consider also the notations:
$$
\sigma^2=\int_{\R^d}K^2(u)\,du\quad\textrm{and}\quad
\eta=\ds{\sum_{k\in\Z^d}E(\varepsilon_{0}\varepsilon_{k}\vert\I)}.
$$
The following proposition (see \cite{Elmachkouri-SISP}) gives the convergence of $Eg_{n}(x)$ to
$g(x)$.
\begin{Prop}\label{rate-biais}
Assume that the assumption $\emph{\textbf{A2)}}$ holds then
$$
\sup_{x\in [0,1]^{d}}\sup_{g\in\C^{1}(B)}\vert
Eg_{n}(x)-g(x)\vert=O\left[h_{n}\right].
$$
\end{Prop}
By proposition 3 in \cite{JD-tcl}, we know that under condition
($\ref{proj}$), the random variable $\eta$ belongs to $L^1$. Our
main result is the following.\\
\vspace{-0.2cm}
\\
\textbf{Main theorem.} {\em If $nh_{n}^{d+1}\to\infty$ and the condition $(\ref{proj})$ holds
then for any $k\in\N^{\ast}$ and any distinct points $x_1,...,x_k$ in $[0,1]^d$, the sequence
$$
(nh_{n})^{d/2}
\left(\begin{array}{c}
       g_{n}(x_1)-Eg_n(x_1)\\
       \vdots\\
       g_{n}(x_k)-Eg_n(x_k)
       \end{array} \right)
\converge{n}{+\infty}{\textrm{$\mathcal{L}$}}
\sigma\,\sqrt{\eta}
\left(\begin{array}{c}
       \tau^{(1)}\\
       \vdots\\
       \tau^{(k)}
       \end{array} \right)\quad\textrm{($\I$-stably)}
$$
where $\sigma^2=\int_{\R^d}K^2(u)\,du$ and $(\tau^{(i)})_{1\leq i\leq k}\sim\mathcal{N}(0,\mathbb{I}_k)$ where $\mathbb{I}_k$ is the identity matrix. Moreover, $(\tau^{(i)})_{1\leq i\leq k}$ is independent of $\eta=\sum_{k\in\Z^d}E(\varepsilon_{0}\varepsilon_{k}\vert\I)$.}\\
\vspace{-0.2cm}
\\
As a consequence of this theorem, we obtain the following result for strongly mixing random fields.\\
\vspace{-0.2cm}
\\
\textbf{Corollary.} {\em Let us consider the following assumption}
\begin{equation}\label{condition-quantile}
\sum_{k\in\Z^d}\int_{0}^{\alpha_{1,\infty}(\vert
k\vert)}\mathcal{Q}_{\varepsilon_{0}}^2(u)\,du<\infty
\end{equation}
{\em where $Q_{\varepsilon_{0}}$ denotes the cadlag inverse of the
function
$H_{\varepsilon_{0}}:t\to\P\left(\vert\varepsilon_{0}\vert>t
\right)$. Then $(\ref{condition-quantile})$ implies $(\ref{proj})$
and also the main theorem.}\\
\vspace{-0.2cm}
\\
\textbf{Remark.} If $\varepsilon_{0}$ is $(2+\delta)$-integrable for some $\delta>0$
then the condition
\begin{equation}\label{condition-sm-strong-mixing}
\sum_{m=1}^{\infty}m^{d-1}\alpha_{1,\infty}^{\delta/(2+\delta)}(m)<\infty
\end{equation}
is more restrictive than condition $(\ref{condition-quantile})$.\\
\\
In order to use the main theorem for establishing confidence intervals, one needs to estimate $\eta$. It is done by the following result established in \cite{JD-tcl}.
\begin{Prop}\label{estimation-eta}
Assume that the condition $(\ref{proj})$ holds. For any $N\in\N^{\ast}$, set $G_N=\{(i,j)\in\Lambda_n\times\Lambda_n\,;\,\vert i-j\vert\leq N\}$. Let $\rho_n$ be a sequence of positive integers satisfying:
$$
\lim_{n\to+\infty}\rho_n=+\infty\quad\textrm{and}\quad\lim_{n\to+\infty}\rho_n^{3d}E(\varepsilon_0^2(1\wedge n^{-d}\varepsilon_0^2)=0
$$
Then
$$
\frac{1}{n^d}\max\left(1,\sum_{(i,j)\in G_{\rho_n}}\varepsilon_i\varepsilon_j\right)\converge{n}{+\infty}{$\P$}\eta.
$$
\end{Prop}
%---------------------------------------------------------------------------------------------------------------------------
\section{Proofs}
\subsection{Proof of the main theorem}
Let $x$ in $[0,1]^{d}$ and $n\geq 1$ be fixed. For any $i$ in
$\Lambda_{n}$, denote
$$
a_{i}(x)=K\left(\frac{x-i/n}{h_{n}}\right)\quad\textrm{and}\quad
b_{i}(x)=\frac{a_{i}(x)}{\sqrt{\sum_{j\in\Lambda_{n}}a_{j}^2(x)}}.
$$
Denote also
$$
v_n(x)=\sqrt{\frac{(nh_{n})^d}{\sum_{i\in\Lambda_{n}}a_{i}(x)}}
\times\sqrt{\frac{\sum_{i\in\Lambda_{n}}a_{i}^2(x)}{\sum_{i\in\Lambda_{n}}a_{i}(x)}}.
$$
Without loss of generality, we consider the case $k=2$ and we refer to $x_1$ and $x_2$ as $x$ and $y$.
Let $\lambda_1$ and $\lambda_2$ be two real numbers such that $\lambda_1^2+\lambda_2^2=1$ and let $x,y\in[0,1]^d$ such that $x\neq y$.
One can notice that
$$
\frac{(nh_{n})^{d/2}}{\sigma}\left[\lambda_1(g_{n}(x)-Eg_{n}(x))+\lambda_2(g_{n}(y)-Eg_{n}(y))\right]=
\sum_{i\in\Lambda_{n}}\tilde{s}_{i}(x,y)\,\varepsilon_{i}
$$
where $\tilde{s}_i(x,y)=(\lambda_1v_n(x)b_i(x)+\lambda_2v_n(y)b_i(y))/\sigma$.
\begin{lemma}\label{convergence-K}
Let $x,y\in[0,1]^d$ be fixed. If $nh_n^{d+1}\to\infty$ then
\begin{equation}\label{convergence-ax-ay}
\lim_{n\to +\infty}\frac{1}{(nh_{n})^d}\sum_{i\in\Lambda_{n}}a_{i}(x)a_i(y)=\delta_{xy}\,\sigma^2
\end{equation}
and
\begin{equation}\label{convergence-ax}
\lim_{n\to +\infty}\frac{1}{(nh_{n})^d}\sum_{i\in\Lambda_{n}}a_{i}(x)=1
\end{equation}
where $\delta_{xy}$ equals $1$ if $x=y$ and $0$ if $x\neq y$.
\end{lemma}
{\em Proof of Lemma \ref{convergence-K}}. In the sequel, we denote $\psi(u)=\frac{1}{h_{n}^d}K\left(\frac{x-u}{h_{n}}\right)K\left(\frac{y-u}{h_{n}}\right)$ and $I_{n}(x,y)=\int_{[0,1]^d}\psi(u)\,du$, we have
$$
I_{n}(x,y)=\sum_{i\in\Lambda_{n}}\int_{R_{i/n}}\psi(u)\,du=\sum_{i\in\Lambda_{n}}n^{-d}\psi(c_{i})
$$
where $R_{i/n}=](i_1-1)/n,i_1/n]\times...\times](i_d-1)/n,i_d/n]$
and $\lambda$ is the Lebesgue measure on $\R^{d}$. Let
$\varphi_{x}(u)=(x-u)/h_{n}$, for any $v$ in $[0,1]^d$, we have
$$
d(K\circ\varphi_{x})(u)(v)=\frac{-1}{h_{n}}\sum_{i=1}^d
v_{i}\sum_{j=1}^d \frac{\partial K}{\partial
u_{j}}(\varphi_{x}(u)).
$$
Using the assumptions on the kernel $K$ and noting that
$$
d\psi(u)=\frac{1}{h_{n}^d}\bigg[d(K\circ\varphi_{x})(u)\times K(\varphi_{y}(u))+d(K\circ\varphi_{y})(u)\times K(\varphi_{x}(u))\bigg]
$$
we derive that there exists $c>0$ such that
$\sup_{u\in[0,1]^d}\|d\psi(u)\|\leq ch_{n}^{-(d+1)}$. So, it
follows that
\begin{align*}
\bigg\vert\frac{1}{(nh_{n})^d}\sum_{i\in\Lambda_{n}}a_{i}(x)a_i(y)-I_{n}(x,y)\bigg\vert
&=\bigg\vert\sum_{i\in\Lambda_{n}}n^{-d}(\psi(i/n)-\psi(c_{i}))\bigg\vert\\
&\leq\sup_{u\in[0,1]^d}\|d\psi(u)\|\,\sum_{i\in\Lambda_{n}}n^{-d}\|i/n-c_{i}\|_{\infty}\\
&=\frac{c}{nh_{n}^{d+1}}\converge{n}{+\infty}{}0.
\end{align*}
Moreover,
$$
I_n(x,y)=\int_{\varphi_{x}([0,1]^d)}K(u)K\left(u+\frac{y-x}{h_n}\right)\,du.
$$
So, we obtain $\lim_{n\to+\infty}I_n(x,y)=\delta_{xy}\,\sigma^2$ and consequently ($\ref{convergence-ax-ay}$) holds. The proof of ($\ref{convergence-ax}$) follows the same lines.
The proof of Lemma $\ref{convergence-K}$ is complete.$\qquad\qquad\Box$\\
\\
Using Lemma \ref{convergence-K} and denoting $\kappa_{xy}^2=(\lambda_1+\lambda_2)^2\delta_{xy}+1-\delta_{xy}$, we derive
$$
\lim_{n\to+\infty}\sum_{i\in\Lambda_n}\tilde{s}_i^2(x,y)=\kappa_{xy}^2=1\quad\textrm{(since $x\neq y$)}.
$$
So, denoting
$$
s_i(x,y)=\frac{\tilde{s}_i(x,y)}{\sqrt{\sum_{j\in\Lambda_n}\tilde{s}_j^2(x,y)}},
$$
it suffices to prove the convergence $\I$-stably of $\sum_{i\in\Lambda_{n}}s_{i}(x,y)\,\varepsilon_{i}$ to
$\sqrt{\eta}\tau_0$ where $\tau_0\sim\mathcal{N}(0,2)$. In fact,
we are going to adapt the proof of the central limit theorem by Dedecker
\cite{JD-tcl}. For any $i$ in $\Z^d$, let us define the tail
$\sigma$-algebra $\F_{i,-\infty}=\cap_{k\in\N^{\ast}}\F_{V_{i}^k}$
(we are going to note $\F_{-\infty}$ in place of $\F_{0,-\infty}$)
and consider the following proposition established in
\cite{JD-tcl}.\\
\vspace{-0.2cm}
\\
\textbf{Proposition} {\em The $\sigma$-algebra $\I$ is included in the $\P$-completion of $\F_{-\infty}$.}\\
\vspace{-0.2cm}
\\
Let $f$ be a one to one map from $[1,N]\cap\N^{\ast}$ to a finite
subset of $\Z^d$ and $(\xi_i)_{i\in\Z^d}$ a real random field. For
all integers $k$ in $[1,N]$, we denote
$$
S_{f(k)}(\xi)=\sum_{i=1}^k \xi_{f(i)}\quad\textrm{and}\quad
S_{f(k)}^{c}(\xi)=\sum_{i=k}^N \xi_{f(i)}
$$
with the convention $S_{f(0)}(\xi)=S_{f(N+1)}^{c}(\xi)=0$. To
describe the set $\Lambda_{n}=\{1,...,n\}^d$, we define the one to
one map $f_{n}$ from $[1,n^d]\cap\N^{\ast}$ to $\Lambda_{n}$ by:
$f_{n}$ is the unique function such that for $1\leq k<l\leq
n^d,\,f(k)<_{lex}f(l)$. From now on, we consider two independent fields
$(\tau^{(1)}_{i})_{i\in\Z^d}$ and $(\tau^{(2)}_{i})_{i\in\Z^d}$ of i.i.d. random variables independent of
$(\varepsilon_{i})_{i\in\Z^d}$ and $\mathcal{I}$ such that $\tau_{0}^{(1)}$ and $\tau_{0}^{(2)}$ have the standard normal
law $\mathcal{N}(0,1)$. We introduce the two sequences of fields $X_{i}=s_{i}(x,y)\varepsilon_{i}$ and
$\gamma_{i}=s_i(x,y)\tau_i\sqrt{\eta}$ where $\tau_i=\tau_{i}^{(1)}+\tau_{i}^{(2)}\sim\mathcal{N}(0,2)$.
Let $h$ be any function from $\R$ to $\R$. For $0\leq k\leq l\leq n^d +1$, we introduce
$h_{k,l}(X)=h(S_{f(k)}(X)+S_{f(l)}^{c}(\gamma))$. With the above
convention we have that $h_{k,n^d+1}(X)=h(S_{f(k)}(X))$ and also
$h_{0,l}(X)=h(S_{f(l)}^{c}(\gamma))$. In the sequel, we will often
write $h_{k,l}$ instead of $h_{k,l}(X)$ and $s_{i}$ instead of $s_{i
}(x,y)$. We denote by
$B_{1}^4(\R)$ the unit ball of $C_{b}^4(\R)$: $h$ belongs to
$B_{1}^4(\R)$ if and only if it belongs to $C^4(\R)$ and satisfies
$\max_{0\leq i\leq 4}\|h^{(i)}\|_{\infty}\leq 1$.\\
\\
{\em 3.1.1 Lindeberg's decomposition}\\
\\
Let $Z$ be a $\I$-measurable random variable bounded by 1. It suffices to prove that for all $h$ in $B_{1}^4(\R)$,
$$
\lim_{n\to+\infty}E\left(Zh(S_{f(n^d)}(X))\right)=E\left(Zh\left((\lambda_1\tau_{0}^{(1)}+\lambda_2\tau_0^{(2)})\sqrt{\eta}\right)\right).
$$
We use Lindeberg's decomposition:
$$
E\left(Z\left[h(S_{f(n^d)}(X))-h\left((\lambda_1\tau_{0}^{(1)}+\lambda_2\tau_0^{(2)})\sqrt{\eta}\right)\right]\right)
=\sum_{k=1}^{n^d}E\left(Z[h_{k,k+1}-h_{k-1,k}]\right).
$$
Now,
$$
h_{k,k+1}-h_{k-1,k}=h_{k,k+1}-h_{k-1,k+1}+h_{k-1,k+1}-h_{k-1,k}.
$$
Applying Taylor's formula we get that:
$$
h_{k,k+1}-h_{k-1,k+1}=X_{f(k)}h_{k-1,k+1}^{'}+\frac{1}{2}X_{f(k)}^{2}h_{k-1,k+1}^{''}+R_{k}
$$
and
$$
h_{k-1,k+1}-h_{k-1,k}=-\gamma_{f(k)}h_{k-1,k+1}^{'}-\frac{1}{2}\gamma_{f(k)}^{2}h_{k-1,k+1}^{''}+r_{k}
$$
where $\vert R_{k}\vert\leq X_{f(k)}^2(1\wedge\vert
X_{f(k)}\vert)$ and $\vert r_{k}\vert\leq\gamma_{f(k)}^2(1\wedge\vert\gamma_{f(k)}\vert)$.
Since $(X,\tau_{i})_{i\neq f(k)}$ is independent of
$\tau_{f(k)}$, it follows that
$$
E\left(Z\gamma_{f(k)}h_{k-1,k+1}^{'}\right)=0\quad\textrm{and}\quad
E\left(Z\gamma_{f(k)}^2h_{k-1,k+1}^{''}\right)=E\left(Zs_{f(k)}^2\eta h_{k-1,k+1}^{''}\right)
$$
Hence, we obtain
\begin{align*}
E\left(Z\left[h(S_{n}(X))-h\left((\lambda_1\tau_{0}^{(1)}+\lambda_2\tau_0^{(2)})\sqrt{\eta}\right)\right]\right)&=
\sum_{k=1}^{n^d}E(ZX_{f(k)}h_{k-1,k+1}^{'})\\
&\quad+\sum_{k=1}^{n^d}E\left(Z\left(X_{f(k)}^2-s_{f(k)}^2\eta\right)\frac{h_{k-1,k+1}^{''}}{2}\right)\\
&\quad+\sum_{k=1}^{n^d}E\left(R_{k}+r_{k}\right).
\end{align*}
Arguing as in Rio \cite{Rio95}, it is proved that $\lim_{n\to+\infty}\sum_{k=1}^{n^d}E\left(\vert R_{k}\vert+\vert
r_{k}\vert\right)=0$. Let us denote $C_N=[-N,N]^d\cap\Z^d$ for any positive integer $N$.
If we define $\eta_{N}=\sum_{k\in C_{N-1}}E\left(\varepsilon_{0}\varepsilon_{k}\vert\I\right)$, the
upper bound $E\vert\eta-\eta_{N}\vert\leq 2\sum_{k\in V_{0}^N} E\vert E\left(\varepsilon_{0}\varepsilon_{k}\vert\I\right)\vert$
holds. Hence according to condition ($\ref{proj}$) and the above
proposition, we derive
$\lim_{N\to+\infty}E\vert\eta-\eta_{N}\vert=0$ and consequently we
have only to show
\begin{equation}\label{equation1}
\lim_{N\to+\infty}\limsup_{n\to+\infty}\sum_{k=1}^{n^d}\left(E(ZX_{f(k)}h_{k-1,k+1}^{'})+E\left(Z\left(X_{f(k)}^2-s_{f(k)}^2\eta_{N}\right)\frac{h_{k-1,k+1}^{''}}{2}\right)\right)=0.
\end{equation}
{\em 3.1.2 First reduction}\\
\\
First, we focus on
$\sum_{k=1}^{n^d}E\left(ZX_{f(k)}h_{k-1,k+1}^{'}\right)$. For all
$N$ in $\N^{\ast}$ and all integer $k$ in $[1,n^d]$, we define
$$
E_{k}^N=f([1,k]\cap\N^{\ast})\cap V_{f(k)}^N\quad\textrm{and}\quad
S_{f(k)}^N(X)=\sum_{i\in E_{k}^N}X_{i}.
$$
For any function $\Psi$ from $\R$ to $\R$, we define
$\Psi_{k-1,l}^N=\Psi(S_{f(k)}^N(X)+S_{f(l)}^c(\gamma))$ (we shall
apply this notation to the successive derivatives of the function
$h$). Our aim is to show that
\begin{equation}\label{equation2}
\lim_{N\to+\infty}\limsup_{n\to+\infty}\sum_{k=1}^{n^d}E\left(Z\left(X_{f(k)}h_{k-1,k+1}^{'}-X_{f(k)}\left(S_{f(k-1)}(X)-S_{f(k)}^N(X)\right)h_{k-1,k+1}^{''}\right)\right)=0.
\end{equation}
First, we use the decomposition
$$
X_{f(k)}h_{k-1,k+1}^{'}=X_{f(k)}h_{k-1,k+1}^{'N}+X_{f(k)}\left(h_{k-1,k+1}^{'}-h_{k-1,k+1}^{'N}\right).
$$
We consider a one to one map $m$ from $[1,\vert
E_{k}^N\vert]\cap\N^{\ast}$ to $E_{k}^N$ and such that $\vert
m(i)-f(k)\vert\leq\vert m(i-1)-f(k)\vert$. This choice of $m$
ensures that $S_{m(i)}(X)$ and $S_{m(i-1)}(X)$ are
$\F_{V_{f(k)}^{\vert m(i)-f(k)\vert}}$-measurable. The fact that
$\gamma$ is independent of $X$ together with proposition 3 in
\cite{JD-tcl} imply that
$$
E\left(ZX_{f(k)}h^{'}\left(S_{f(k+1)}^c(\gamma)\right)\right)=E\left(h^{'}\left(S_{f(k+1)}^c(\gamma)\right)\right)
E\left(ZE\left(X_{f(k)}\vert\F_{-\infty}\right)\right)=0.
$$
Therefore $\vert E\left(ZX_{f(k)}h_{k-1,k+1}^{'N}\right)\vert$
equals
$$
\bigg\vert\sum_{i=1}^{\vert
E_{k}^N\vert}E\left(ZX_{f(k)}\bigg[h^{'}\left(S_{m(i)}(X)+S_{f(k+1)}^c(\gamma)\right)-h^{'}\left(S_{m(i-1)}(X)+S_{f(k+1)}^c(\gamma)\right)\bigg]\right)\bigg\vert.
$$
Since $S_{m(i)}(X)$ and $S_{m(i-1)}(X)$ are $\F_{V_{f(k)}^{\vert
m(i)-f(k)\vert}}$-measurable, we can take the conditional
expectation of $X_{f(k)}$ with respect to $\F_{V_{f(k)}^{\vert
m(i)-f(k)\vert}}$ in the right hand side of the above equation. On
the other hand the function $h^{'}$ is $1$-Lipschitz, hence
$$
\vert
h^{'}\left(S_{m(i)}(X)+S_{f(k+1)}^c(\gamma)\right)-h^{'}\left(S_{m(i-1)}(X)+S_{f(k+1)}^c(\gamma)\right)\vert\leq\vert
X_{m(i)}\vert.
$$
Consequently, the term
$$
\bigg\vert
E\left(ZX_{f(k)}\bigg[h^{'}\left(S_{m(i)}(X)+S_{f(k+1)}^c(\gamma)\right)-h^{'}\left(S_{m(i-1)}(X)+S_{f(k+1)}^c(\gamma)\right)\bigg]\right)\bigg\vert
$$
is bounded by
$$
E\vert X_{m(i)}E_{\vert m(i)-f(k)\vert}\left(X_{f(k)}\right)\vert
$$
and
$$
\vert
E\left(ZX_{f(k)}h_{k-1,k+1}^{'N}\right)\vert\leq\sum_{i=1}^{\vert
E_{k}^N\vert}E\vert X_{m(i)}E_{\vert
m(i)-f(k)\vert}(X_{f(k)})\vert.
$$
Hence,
\begin{align*}
\bigg\vert\sum_{k=1}^{n^d}E\left(ZX_{f(k)}h_{k-1,k+1}^{'N}\right)\bigg\vert
&\leq\sum_{k=1}^{n^d}\vert s_{f(k)}\vert\sum_{i=1}^{\vert
E_{k}^N\vert}\vert s_{m(i)}\vert E\vert\varepsilon_{m(i)}E_{\vert
m(i)-f(k)\vert}(\varepsilon_{f(k)})\vert\\
&\leq A\sum_{j\in V_0^N}\|\varepsilon_jE_{\vert j\vert}(\varepsilon_0)\|_1<+\infty\quad (A\in\R_+^{\ast})
\end{align*}
where (by Lemma \ref{convergence-K}) we used the fact that
\begin{equation}\label{estimation-s-1}
\sup_{i\in\Lambda_n}\vert s_i\vert =O\left(\frac{1}{(nh_n)^{d/2}}\right)
\end{equation}
and
\begin{equation}\label{estimation-s-2}
\sum_{i\in\Lambda_n}\vert s_i\vert=O\left((nh_n)^{d/2}\right).
\end{equation}
Since ($\ref{proj}$) is satisfied, this last term is as small as we
wish by choosing $N$ large enough. Applying again Taylor's
formula, it remains to consider
$$
X_{f(k)}(h_{k-1,k+1}^{'}-h_{k-1,k+1}^{'N})=X_{f(k)}(S_{f(k-1)}(X)-S_{f(k)}^N(X))h_{k-1,k+1}^{''}+R_{k}^{'},
$$
where $\vert R_{k}^{'}\vert\leq 2\vert
X_{f(k)}(S_{f(k-1)}(X)-S_{f(k)}^N(X))(1\wedge\vert
S_{f(k-1)}(X)-S_{f(k)}^N(X)\vert)\vert$. It follows that
\begin{align*}
\sum_{k=1}^{n^d}E\vert R_{k}^{'}\vert &\leq 2A\,E\left(\vert\varepsilon_{0}\vert\left(\sum_{i\in\Lambda_{N}}\vert\varepsilon_{i}\vert\right)
\left(1\wedge\sum_{i\in\Lambda_{N}}\vert s_{i}\vert\vert\varepsilon_{i}\vert\right)\right)\quad (A\in\R_+^{\ast}).
\end{align*}
Keeping in mind that $s_{i}\to 0$ as $n\to\infty$ and applying the dominated convergence theorem,
this last term converges to zero as $n$ tends to infinity and ($\ref{equation2}$) follows.\\
\\
{\em 3.1.3 The second order terms}\\
\\
It remains to control
\begin{equation}\label{W}
W_{1}=E\left(Z\sum_{k=1}^{n^d}h_{k-1,k+1}^{''}\left(\frac{X_{f(k)}^2}{2}+X_{f(k)}\left(S_{f(k-1)}(X)-S_{f(k)}^N(X)\right)-\frac{s_{f(k)}^2\eta_{N}}{2}\right)\right).
\end{equation}
We consider the following sets:
$$
\Lambda_{n}^N=\{i\in\Lambda_{n}\,;\,d({i},\partial\Lambda_{n})\geq
N\}\quad\textrm{and}\quad I_{n}^N=\{1\leq i\leq n^d\,;\,f(i)\in
\Lambda_{n}^N\},
$$
and the function $\Psi$ from $\R^{\Z^d}$ to $\R$ such that
$$
\Psi(\varepsilon)=\varepsilon_{0}^2+\sum_{i\in V_{0}^1\cap
C_{N-1}}2\varepsilon_{0}\varepsilon_{i}.
$$
For $k$ in $[1,n^d]$, we set $D_{k}^N=\eta_{N}-\Psi\circ T^{f(k)}(\varepsilon)$. By definition of $\Psi$ and of the set
$I_{n}^N$, we have for any $k$ in $I_{n}^N$
$$
\Psi\circ
T^{f(k)}(\varepsilon)=\varepsilon_{f(k)}^2+2\varepsilon_{f(k)}(S_{f(k-1)}(\varepsilon)-S_{f(k)}^N(\varepsilon)).
$$
Therefore for $k$ in $I_{n}^N$
$$
s_{f(k)}^2D_{k}^N=s_{f(k)}^2\eta_{N}-X_{f(k)}^2-2X_{f(k)}(S_{f(k-1)}(X)-S_{f(k)}^N(X)).
$$
Since $\lim_{n\to+\infty}n^{-d}\vert I_{n}^N\vert=1$, it remains
to prove that
\begin{equation}\label{equation3}
\lim_{N\to+\infty}\limsup_{n\to+\infty}E\left(Z\sum_{k=1}^{n^d}s_{f(k)}^2h_{k-1,k+1}^{''}D_{k}^N\right)=0.
\end{equation}
{\em 3.1.4 Conditional expectation with respect to the tail
$\sigma$-algebra}\\
\\
Now, we are going to replace $D_{k}^N$ by
$E\left(D_{k}^N\vert\F_{f(k),-\infty}\right)$. We introduce the
expression
$$
H_{n}^N=\sum_{k=1}^{n^d}E\left(s_{f(k)}^2Zh_{k-1,k+1}^{''}[\Psi\circ
T^{f(k)}(\varepsilon)-E(\Psi\circ
T^{f(k)}(\varepsilon)\vert\F_{f(k),-\infty})]\right).
$$
For sake of brevity, we have written $h_{k-1,k+1}^{''}$ instead of
$h_{k-1,k+1}^{''}(X)$. Using the stationarity of the field we get
that
$$
H_{n}^N=\sum_{k=1}^{n^d}E\left(s_{f(k)}^2Z(h_{k-1,k+1}^{''}\circ
T^{-f(k)})(X)[\Psi(\varepsilon)-E(\Psi(\varepsilon)\vert\F_{-\infty})]\right).
$$
For any positive integer $p$, we decompose $H_{n}^N$ in two parts
$$
H_{n}^N=\sum_{k=1}^{n^d}J_{k}^1(p)+\sum_{k=1}^{n^d}J_{k}^2(p),
$$
where
$$
J_{k}^1(p)=E\left(s_{f(k)}^2Z(h_{k-1,k+1}^{''p}\circ
T^{-f(k)})[\Psi(\varepsilon)-E(\Psi(\varepsilon)\vert\F_{-\infty})]\right)
$$
and $J_{k}^2(p)$ equals to
$$
E\left(s_{f(k)}^2Z[h_{k-1,k+1}^{''}\circ
T^{-f(k)}-h_{k-1,k+1}^{''p}\circ
T^{-f(k)}](X)[\Psi(\varepsilon)-E(\Psi(\varepsilon)\vert\F_{-\infty})]\right).
$$
From the definition of $h_{k-1,k+1}^{''p}$, we infer that the
variable $h_{k-1,k+1}^{''p}\circ T^{-f(k)}(X)$ is
$\F_{V_{0}^p}$-measurable. Therefore, we can take the conditional
expectation of
$\Psi(\varepsilon)-E(\Psi(\varepsilon)\vert\F_{-\infty})$ with
respect to $\F_{V_{0}^p}$ in the expression of $J_{k}^1(p)$. Now,
the backward martingale limit theorem implies that
$$
\lim_{p\to+\infty}E\vert
E(\Psi(\varepsilon)\vert\F_{V_{0}^p})-E(\Psi(\varepsilon)\vert\F_{-\infty})\vert=0
$$
and consequently
$$
\lim_{p\to+\infty}\limsup_{n\to+\infty}\bigg\vert\sum_{k=1}^{n^d}J_{k}^1(p)\bigg\vert=0.
$$
On the other hand
$$
\bigg\vert\sum_{k=1}^{n^d}J_{k}^2(p)\bigg\vert\leq
E\bigg[\left(2\wedge\sum_{\vert
i\vert<p} s_{f(i)}^2\vert\varepsilon_{i}\vert\right)\vert
\Psi(\varepsilon)-E(\Psi(\varepsilon)\vert\F_{-\infty})\vert\bigg].
$$
Hence, applying the dominated convergence theorem, we conclude that $H_{n}^N$ tends to zero as $n$ tends to infinity. It remains to consider
$$
W_{2}=E\left(Z\sum_{k=1}^{n^d}h_{k-1,k+1}^{''}s_{f(k)}^2E(D_{k}^N\vert\F_{f(k),-\infty})\right).
$$
{\em 3.1.5 Truncation}\\
\\
For any integer $k$ in $[1,n^d]$ and any $M$ in $\R^{+}$ we introduce the two sets
$$
B_{k}^N(M)=E(D_{k}^N\vert\F_{f(k),-\infty})\ind{\vert\eta_{N}-E(\Psi\circ T^{f(k)}(\varepsilon)\vert\F_{f(k),-\infty})\vert\leq M}
$$
and
$$
\overline{B}_{k}^N(M)=E(D_{k}^N\vert\F_{f(k),-\infty})-B_{k}^N(M).
$$
The stationarity of the field ensures that
$E\vert\overline{B}_{k}^N(M)\vert=E\vert\overline{B}_{1}^N(M)\vert$
for any $k$ in $[1,n^d]$. Now, applying the dominated convergence
theorem, we have
$\lim_{M\to+\infty}E\vert\overline{B}_{1}^N(M)\vert=0$. It follows
that
$$
\lim_{M\to+\infty}\sum_{k=1}^{n^d}E\left(h_{k-1,k+1}^{''}s_{f(k)}^2\overline{B}_{k}^N(M)\right)=0.
$$
Therefore instead of $W_{2}$ it remains to consider
$$
W_{3}=E\left(Z\sum_{k=1}^{n^d}h_{k-1,k+1}^{''}s_{f(k)}^2B_{k}^N(M)\right).
$$
{\em 3.1.6 An ergodic lemma}\\
\\
The next result is the central point of the proof.
\begin{lemma}\label{lln}
For all $M$ in $\R^{+}$, we introduce
$$
\beta_{N}(M)=E\left([\eta_{N}-E\left(\Psi(\varepsilon)\vert\F_{-\infty}\right)]\ind{\vert\eta_{N}-E(\Psi(\varepsilon)\vert\F_{-\infty})\vert\leq
M}\big\vert\I\right).
$$
Then
$$
\lim_{M\to+\infty}\beta_{N}(M)=0\quad\textrm{a.s.}\quad\textrm{and}\quad
\lim_{n\to+\infty}E\bigg\vert\beta_{N}(M)-\sum_{k=1}^{n^d}s_{f(k)}^2B_{k}^N(M)\bigg\vert=0.
$$
\end{lemma}
{\em Proof of Lemma $\ref{lln}$.} Let
$$
u(\varepsilon)=[\eta_N-E\left(\Psi(\varepsilon)\vert\F_{-\infty}\right)]\ind{\vert\eta_{N}-E(\Psi(\varepsilon)\vert\F_{-\infty})\vert\leq
M}.
$$
Using the function $u$, we write
$\beta_N(M)=E(u(\varepsilon)\vert\I)$. The fact that $\beta_N(M)$
tends to zero as $M$ tends to infinity follows from the dominated
convergence theorem. In fact $\lim_{M\to\infty}u(\varepsilon)=\eta_{N}-E(\Psi(\varepsilon)\vert\F_{-\infty})$ and $u(\varepsilon)$ is bounded by
$\vert\eta_N-E(\Psi(\varepsilon)\vert\F_{-\infty})\vert$ which
belongs to $L^1$. This implies that
$$
\lim_{M\to\infty}\beta_N(M)=E\left(\eta_N-E(\Psi(\varepsilon)\vert\F_{-\infty})\,\big\vert\I\right)\quad\textrm{a.s.}
$$
Since $\I$ is included in the $\P$-completion of $\F_{-\infty}$
(see the above proposition) and keeping in mind
that $\eta_N$ is $\I$-measurable, it follows that
$$
\lim_{M\to\infty}\beta_N(M)=\eta_{N}-E(\Psi(\varepsilon)\vert\I)\quad\textrm{a.s.}
$$
By stationarity of the random field, we know that
$E(\varepsilon_{0}\varepsilon_{k}\vert\I)=E(\varepsilon_{0}\varepsilon_{-k}\vert\I)$
which implies that $E(\Psi(\varepsilon)\vert\I)=\eta_{N}$ and the result follows.\\
We are going to prove the second point of Lemma $\ref{lln}$. Noting that $B_{k}(M)=u\circ T^{f(k)}(\varepsilon)$, we have
$$
\sum_{k=1}^{n^d}s_{f(k)}^2 B_{k}^N(M)=\sum_{i\in\Lambda_{n}}s_{i}^2\,u\circ T^{i}(\varepsilon).
$$
Finally, the proof of lemma $\ref{lln}$ is completed by the following lemma which the proof is left to the reader.
\begin{lemma}\label{VonNeuman}
$$
\lim_{n\to\infty}\bigg\|\sum_{i\in\Lambda_{n}}s_{i}^2\,u\circ
T^{i}(\varepsilon)-E(u(\varepsilon)\vert\I)\bigg\|_{2}=0.
$$
\end{lemma}
As a direct application of lemma \ref{lln}, we see that
$$
\bigg\vert
E\left(Z\sum_{k=1}^{n^d}h_{k-1,k+1}^{''}s_{f(k)}^2\beta_{N}(M)\right)\bigg\vert\leq
E\vert\beta_{N}(M)\vert
$$
is as small as we wish by choosing $M$ large enough. So instead of
$W_{3}$ we consider
$$
W_{4}=E\left(Z\sum_{k=1}^{n^d}h_{k-1,k+1}^{''}s_{f(k)}^2[B_{k}^N(M)-\beta_{N}(M)]\right).
$$
{\em 3.1.7 Abel transformation}\\
\\
In order to control $W_{4}$, we use the Abel transformation:
\begin{align*}
W_{4}&=E\bigg[\sum_{k=1}^{n^d}\left(\sum_{i=1}^k s_{f(i)}^2[B_{i}^N(M)-\beta_{N}(M)]\right)Z(h_{k-1,k+1}^{''}-h_{k,k+2}^{''})\bigg]\\
&\qquad\qquad
+E\left(Zh_{n^d,n^d+2}^{''}\sum_{k=1}^{n^d}s_{f(k)}^2[B_{k}^N(M)-\beta_{N}(M)]\right).
\end{align*}
Now
$$
\bigg\vert
E\left(Zh_{n^d,n^d+2}^{''}\sum_{k=1}^{n^d}s_{f(k)}^2[B_{k}^N(M)-\beta_{N}(M)]\right)\bigg\vert\leq
E\bigg\vert\beta_{N}(M)-\sum_{k=1}^{n^d}s_{f(k)}^2B_{k}^N(M)\bigg\vert.
$$
Then applying lemma $\ref{lln}$, we obtain
$$
\lim_{n\to+\infty}\bigg\vert
E\left(Zh_{n^d,n^d+2}^{''}\sum_{k=1}^{n^d}s_{f(k)}^2[B_{k}^N(M)-\beta_{N}(M)]\right)\bigg\vert=0.
$$
Therefore it remains to prove that for any positive integer $N$ and any positive real $M$,
$$
\lim_{n\to+\infty} E\bigg[\sum_{k=1}^{n^d}\left(\sum_{i=1}^k
s_{f(i)}^2[B_{i}^N(M)-\beta_{N}(M)]\right)Z(h_{k-1,k+1}^{''}-h_{k,k+2}^{''})\bigg]=0.
$$
{\em 3.1.8 Last reductions}\\
\\
We are going to finish the proof. We use the same decomposition as
before:
$$
h_{k,k+2}^{''}-h_{k-1,k+1}^{''}=h_{k,k+2}^{''}-h_{k,k+1}^{''}+h_{k,k+1}^{''}-h_{k-1,k+1}^{''}.
$$
Applying Taylor's formula, we have $h_{k,k+2}^{''}-h_{k,k+1}^{''}=-\gamma_{f(k+1)}h_{k,k+2}^{'''}+t_{k}$ and $h_{k,k+1}^{''}-h_{k-1,k+1}^{''}=X_{f(k)}h_{k-1,k+1}^{'''}+T_{k}$ where $\vert t_{k}\vert\leq\gamma_{f(k+1)}^2$ and $\vert T_{k}\vert\leq X_{f(k)}^2$. To examine the remainder terms, we consider:
$$
E\left(\sum_{k=1}^{n^d}s_{f(k)}^2\left(\sum_{i=1}^k
s_{f(i)}^2[B_{i}^N(M)-\beta_{N}(M)]\right)Z\varepsilon_{f(k)}^2\right).
$$
The definition of $B_{i}^N(M)$ and of $\beta_{N}(M)$ enables us to write for all integer $k$ in $[1,n^d]$,
$$
\sum_{i=1}^k s_{f(i)}^2\vert B_{i}^N(M)-\beta_{N}(M)\vert\leq
2M.
$$
Therefore
$$
E\bigg\vert\sum_{k=1}^{n^d}\left(\sum_{i=1}^k
s_{f(i)}^2[B_{i}^N(M)-\beta_{n}(M)]\right)s_{f(k)}^2
Z\varepsilon_{f(k)}^2\ind{\vert\varepsilon_{f(k)}\vert>K}\bigg\vert\leq
2ME\left(\varepsilon_{0}^2\ind{\vert\varepsilon_0\vert>K}\right)
$$
and applying the dominated convergence theorem this last term is as small as we wish by choosing $K$ large enough. Now, for all $K$
in $\R^{+}$, Lemma $\ref{lln}$ ensures that
$$
\lim_{n\to+\infty}E\left(\sum_{k=1}^{n^d}s_{f(k)}^2\left(\sum_{i=1}^k
s_{f(i)}^2[B_{i}^N(M)-\beta_{N}(M)]\right)Z\varepsilon_{f(k)}^2\ind{\vert\varepsilon_{f(k)}\vert\leq K}\right)=0.
$$
So, we have proved that
$$
\lim_{n\to+\infty} E\left(\sum_{k=1}^{n^d}\left(\sum_{i=1}^k
s_{f(i)}^2[B_{i}^N(M)-\beta_{N}(M)]\right)ZT_{k}\right)=0.
$$
In the same way, we obtain that
$$
\lim_{n\to+\infty} E\left(\sum_{k=1}^{n^d}\left(\sum_{i=1}^k
s_{f(i)}^2[B_{i}^N(M)-\beta_{N}(M)]\right)Zt_{k}\right)=0.
$$
Moreover since $(\varepsilon,(\tau_{i})_{i\neq f(k+1)})$ is
independent of $\tau_{f(k+1)}$ we have
$$
E\left(\sum_{i=1}^k s_{f(i)}^2[B_{i}^N(M)-\beta_{N}(M)]\gamma_{f(k+1)}Zh_{k,k+2}^{'''}\right)=0.
$$
Finally, it remains to consider
$$
W_{5}=E\bigg[\sum_{k=1}^{n^d}\left(\sum_{i=1}^k s_{f(i)}^2
[B_{i}^N(M)-\beta_{N}(M)]\right)ZX_{f(k)}h_{k-1,k+1}^{'''}\bigg].
$$
Let $p$ be a fixed positive integer. Since $h^{'''}$ is
$1$-Lipschitz, we have the upper bound $\vert
h_{k-1,k+1}^{'''}-h_{k-1,k+1}^{'''p}\vert\leq\vert
S_{f(k-1)}(X)-S_{f(k)}^p(X)\vert$. Now, we can apply the same
truncation argument as before: first we choose the level of our
truncation by applying the dominated convergence theorem and then
we use Lemma $\ref{lln}$. So, it follows that
$$
\lim_{n\to+\infty}E\bigg[\sum_{k=1}^{n^d}\left(\sum_{i=1}^k
s_{f(i)}^2[B_{i}^N(M)-\beta_{N}(M)]\right)ZX_{f(k)}(h_{k-1,k+1}^{'''}-h_{k-1,k+1}^{'''p})\bigg]=0.
$$
Therefore, to prove our theorem it is enough to show that
\begin{equation}\label{equation4}
\lim_{p\to+\infty}\limsup_{n\to+\infty}E\bigg[\sum_{k=1}^{n^d}\left(\sum_{i=1}^k
s_{f(i)}^2[B_{i}^N(M)-\beta_{N}(M)]\right)ZX_{f(k)}h_{k-1,k+1}^{'''p}\bigg]=0.
\end{equation}
We consider a one to one map $m$ from $[1,\vert
E_{k}^p\vert]\cap\N^{\ast}$ to $E_{k}^p$ and such that $\vert
m(i)-f(k)\vert\leq\vert m(i-1)-f(k)\vert$. Now, we use the same
argument as before:
\begin{align*}
h_{k-1,k+1}^{'''p}-h^{'''}(S_{f(k)}^c(\gamma))
&=\sum_{i=1}^{\vert
E_{k}^p\vert}h^{'''}(S_{m(i)}(X)+S_{f(k)}^c(\gamma))-h^{'''}(S_{m(i-1)}(X)+S_{f(k)}^c(\gamma))\\
&\leq\sum_{i=1}^{\vert
E_{k}^p\vert}\vert X_{m(i)}\vert.
\end{align*}
Here recall that $B_{i}^N(M)$ is $\F_{f(i),-\infty}$-measurable
and $\beta_{N}(M)$ is $\I$-measurable. We have
$E(\varepsilon_{f(k)}\vert\I)=0,\,E(\varepsilon_{f(k)}\vert\F_{f(k),-\infty})=0$
and $E(\varepsilon_{f(k)}\vert\F_{f(i),-\infty})=0$ for any
positive integer $i$ such that $i<k$. Consequently, for any
positive integer $i$ such that $i\leq k$, we have
$$
E\left(s_{f(i)}^2[B_{i}^N(M)-\beta_{N}(M)]Zs_{f(k)}\varepsilon_{f(k)}h^{'''}(S_{f(k)}^c(\gamma))\right)=0.
$$
Therefore using the conditional expectation, we find
\begin{align*}
&E\bigg[\sum_{k=1}^{n^d}\left(\sum_{i=1}^k
s_{f(i)}^2[B_{i}^N(M)-\beta_{N}(M)]\right)ZX_{f(k)}h^{'''p}_{k-1,k+1}\bigg]\\
&\leq 2M\sum_{k=1}^{n^d}\vert s_{f(k)}\vert\sum_{i=1}^{\vert
E_{k}^p\vert}\vert s_{m(i)}\vert E\vert\varepsilon_{m(i)}E_{\vert
m(i)-f(k)\vert}(\varepsilon_{f(k)})\vert\\
&\quad=2M\sum_{k=1}^{n^d}\vert s_{f(k)}\vert\sum_{j\in
V_{0}^p}\vert s_{j+f(k)}\vert E\vert\varepsilon_{j}E_{\vert
j\vert}(\varepsilon_{0})\vert\\
&\qquad\leq 2AM\sum_{j\in V_{0}^p}E\vert\varepsilon_{j}E_{\vert
j\vert}(\varepsilon_{0})\vert\quad (A\in\R_+^{\ast})\quad\textrm{by $(\ref{estimation-s-1})$ and $(\ref{estimation-s-2})$}.
\end{align*}
Since ($\ref{proj}$) is realised the last term is as small as we
wish by choosing $p$ large enough, hence $W_{4}$ is handled.
Finally, the main theorem is proved.$\qquad\qquad\Box$
\subsection{Proof of the corollary}
As observed in \cite{JD-tcl}, the proof of the corollary is a
direct consequence of Theorem 1.1 in Rio \cite{Rio93}. In fact, for any $k$ in $V_{0}^{1}$, we
have
\begin{align*}
E\vert\varepsilon_k E_{\vert k\vert}(\varepsilon_{0})\vert
&\leq 4\int_{0}^{\alpha_{1,\infty}(\vert
k\vert)}Q_{\varepsilon_{0}}^2(u)\,du.
\end{align*}
The proof of the corollary is complete.$\qquad\qquad\Box$

\section{Application}
The direct consequence of our result is that it allows the
construction of statistical tests able to quantify the estimation
error. For this purpose, we show the construction of such a test that
can be used in image denoising~\cite{Gonzalez-Woods,Li,Winkler}. In
the context given by the model~(\ref{model}), let us consider the
following situation~: a true image $g$ is affected by a correlated
additive noise $\epsilon$, that gives $Y$ for the observed image.

For the original function two images were considered. The first one is
a simulated image, a two-dimensional sinusoide, whereas the second one
is the very well known Lena image. The first image since it represents
a continuous function, perfectly matches the hypothesis of our
results. The second one represents a piece-wise continuous function,
so the hypothesis of our result are not completely verified, still
this is a much more realistics situation.

These images are gray levels images with pixels values in the interval
$[0,255]$. The size of an image is $256 \times 256$ pixels. The
correlated noise we consider is a Gaussian field $(\varepsilon_{k})_{k
\in \Z^2}$ built using an exponential covariance function
$$ C(k) =E(\varepsilon_{0}\varepsilon_{k}) = \text{Cst} \times
\exp\{-\frac{|k|}{a}\}.
$$

The choice of such random field ensures the
validity of the projective
criterion~(\ref{proj}) $\quad$ (see~\cite{Doukhan}, p.59, Corollary 2). There
exist several methods for simulating such a random field, here we
have opted for the spectral method~\cite{Lantuejoul}. In order to
obtain an important visual effect of how the noise affects the
original image $\text{Cst}$ was set to $200$ and $a = 1$. The
noisy image is obtained by adding pixel by pixel the original image to
the simulated noise. The estimator of the original image is computed
using the Epanechnikov kernel
$$ K(x) = \frac{3}{8} (1- |x|^{2})\mathbb{I}_{\{|x| \leq
1\}},\,\,x=(x_{1},x_{2}) \in \R^2.
$$

In order to compute the expectation of the estimated function, several
realisation of the noisy image are needed. Here we have considered
$50$ such images, constructed by adding the original Lena image with a
noise realisation. Using~(\ref{def-g_n}), for each noisy image, an
estimate $g_n$ of the original function $g$ was computed using the
kernel $K$ defined above. The expectation $E ( g_n )$ is computed
by just taking the pixel by pixel arithmetical means corresponding to
the images previously restored.

Clearly, it is now possible to estimate the difference $g_n
- E(g_n)$. Following our theoretical result, the normalised square of this
difference follows a $\chi^{2}$ distribution with one
degree of freedom. Since this quantity is observable, $p$-values
pixel by pixel can be computed.

The obtained results for the synthetic and real image restoration are
shown in Figure~\ref{results_synth} and \ref{results},
respectively. In both situation, it can be noticed that in the
``dirty'' pictures, spots are formed, due to the noise
correlation. The expectations of the estimated original images exhibit
almost no such spots. Furthermore, the visual quality of the restored
images is close to the originals. A more quantitative evaluation of
this result is given by the image of $p$-values of the proposed
statistical test given in. The light-coloured pixels represent
$p$-values close to $1$, whereas the dark-coloured pixels indicate
values close to $0$. For the real image case, we have counted $83\%$
of the pixels for which we have obtained a $p$-value higher than
$0.01$. This ratio is quite a reliable indicator concerning the
restored image. Together with the visual analysis of the results, it
provides a detailed description of the obtained result. We conclude
that, under these considerations, the theoretical results developed in
this paper may be used as a basis for the development of practical
tools in image analysis.

\begin{figure}[!htbp]
\begin{center}
\begin{tabular}{cc}
a)\epsfxsize=5cm \epsffile{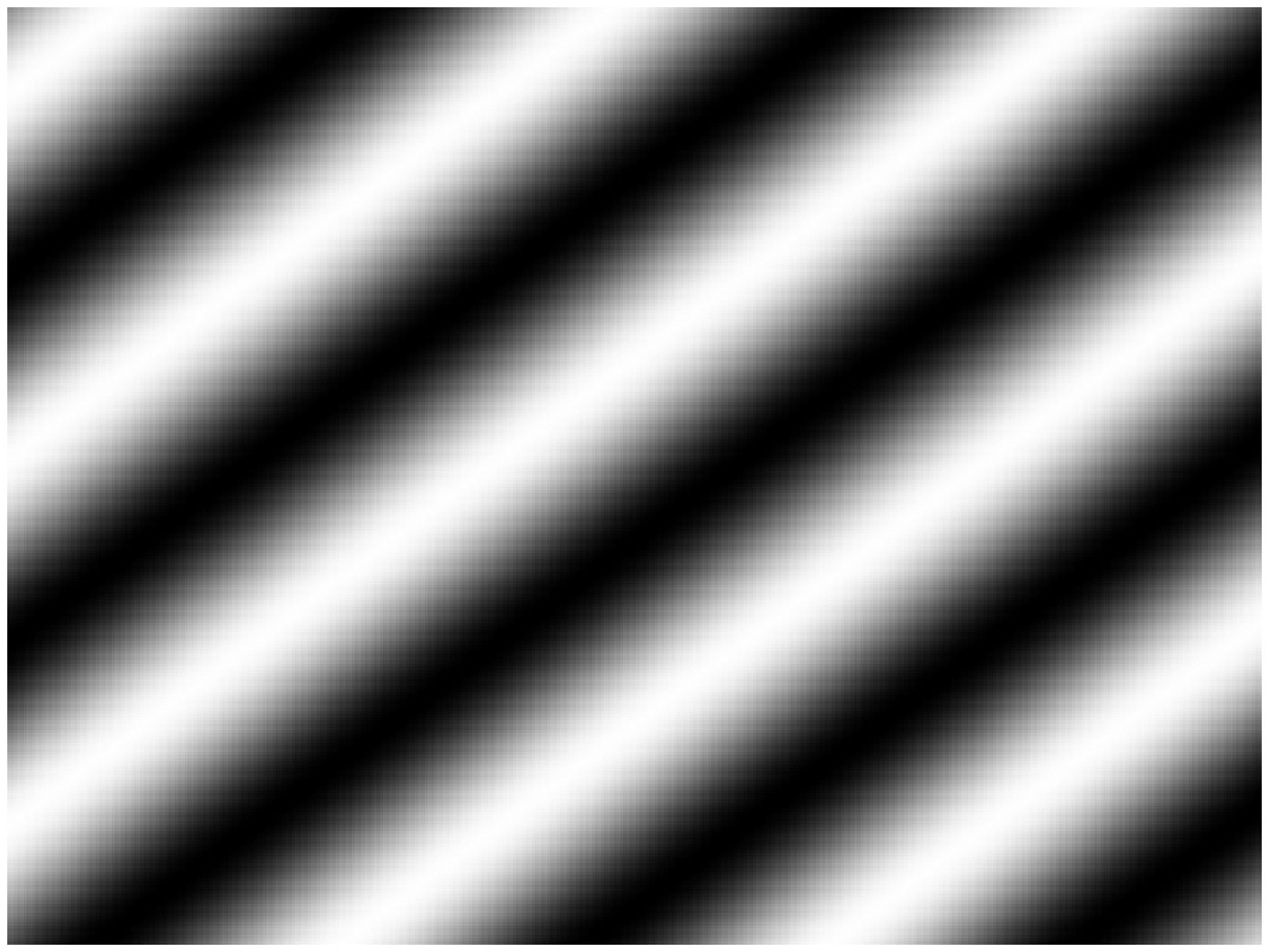} &
b)\epsfxsize=5cm \epsffile{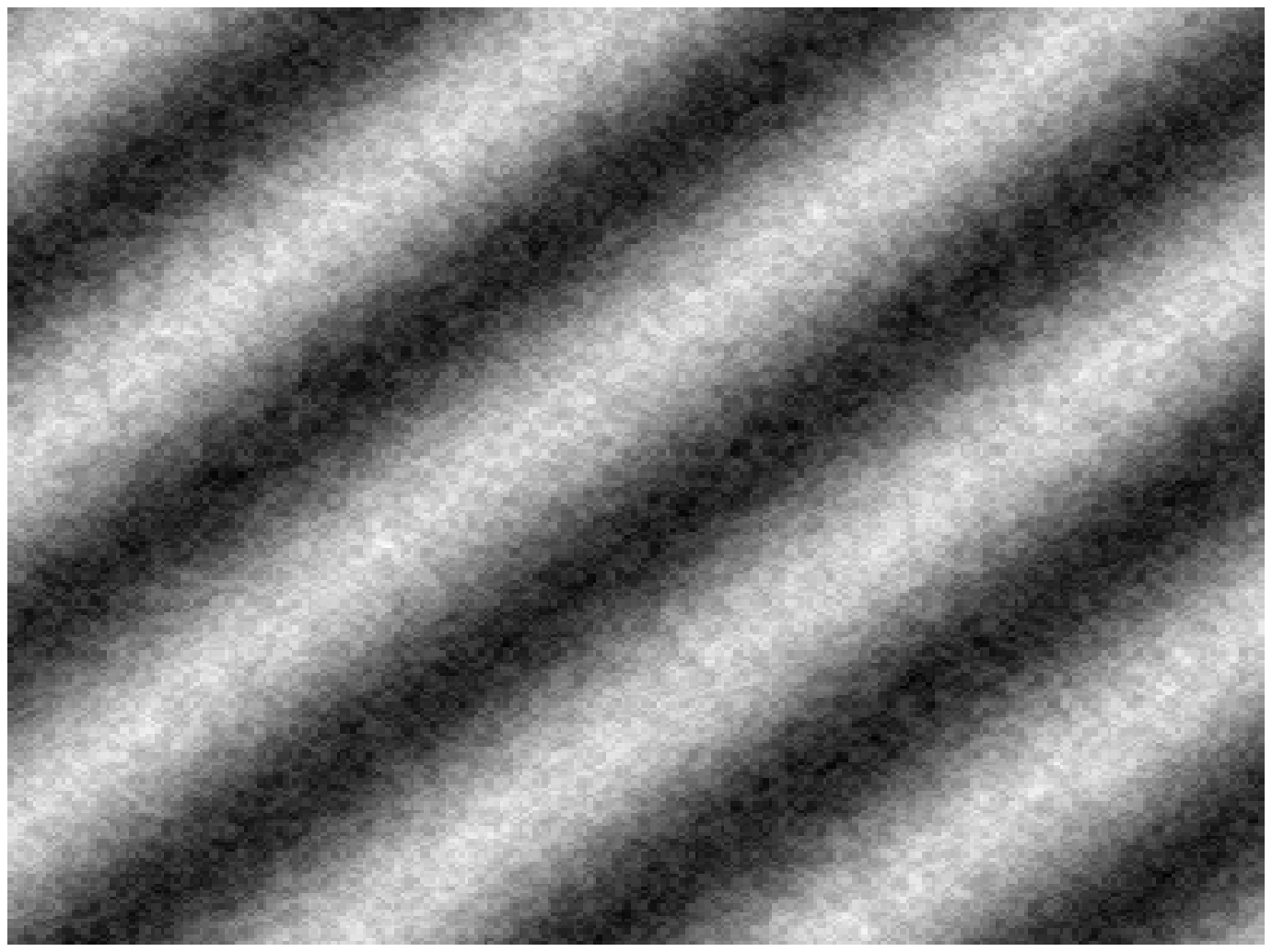}\\
c)\epsfxsize=5cm \epsffile{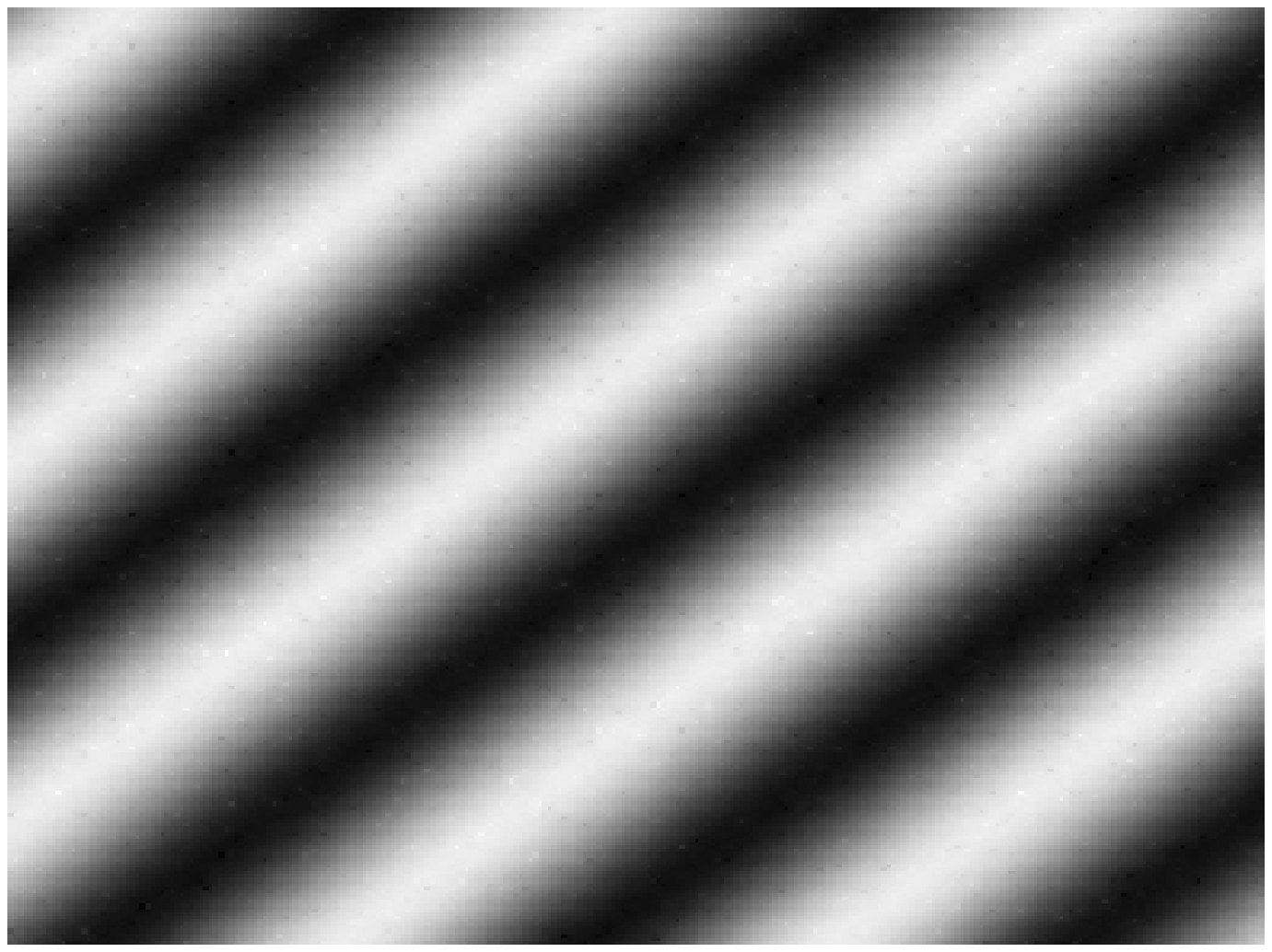} &
d)\epsfxsize=5cm \epsffile{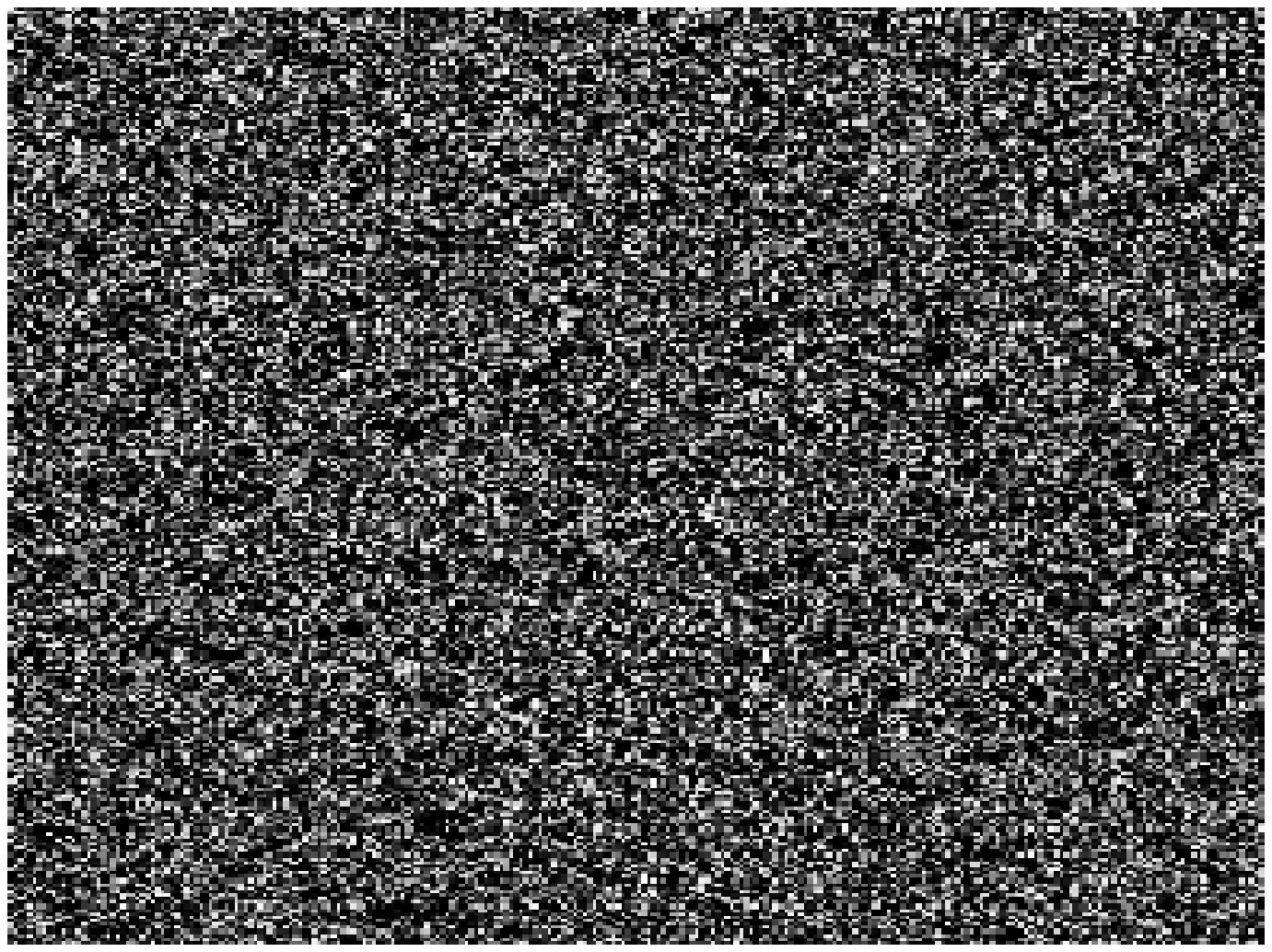}\\
\end{tabular}
\end{center}
\caption{Results of the image restoration procedure~: a) original sinusoide
image, b) realisation of a noisy image, c) expectation of the restored
images, d) obtained $p-$values as a gray level image (white pixels
represent values close to $1$, whereas black pixels indicate values
close to $0$).}
\label{results_synth}
\end{figure}

\begin{figure}[!htbp]
\begin{center}
\begin{tabular}{cc}
a)\epsfxsize=5cm \epsffile{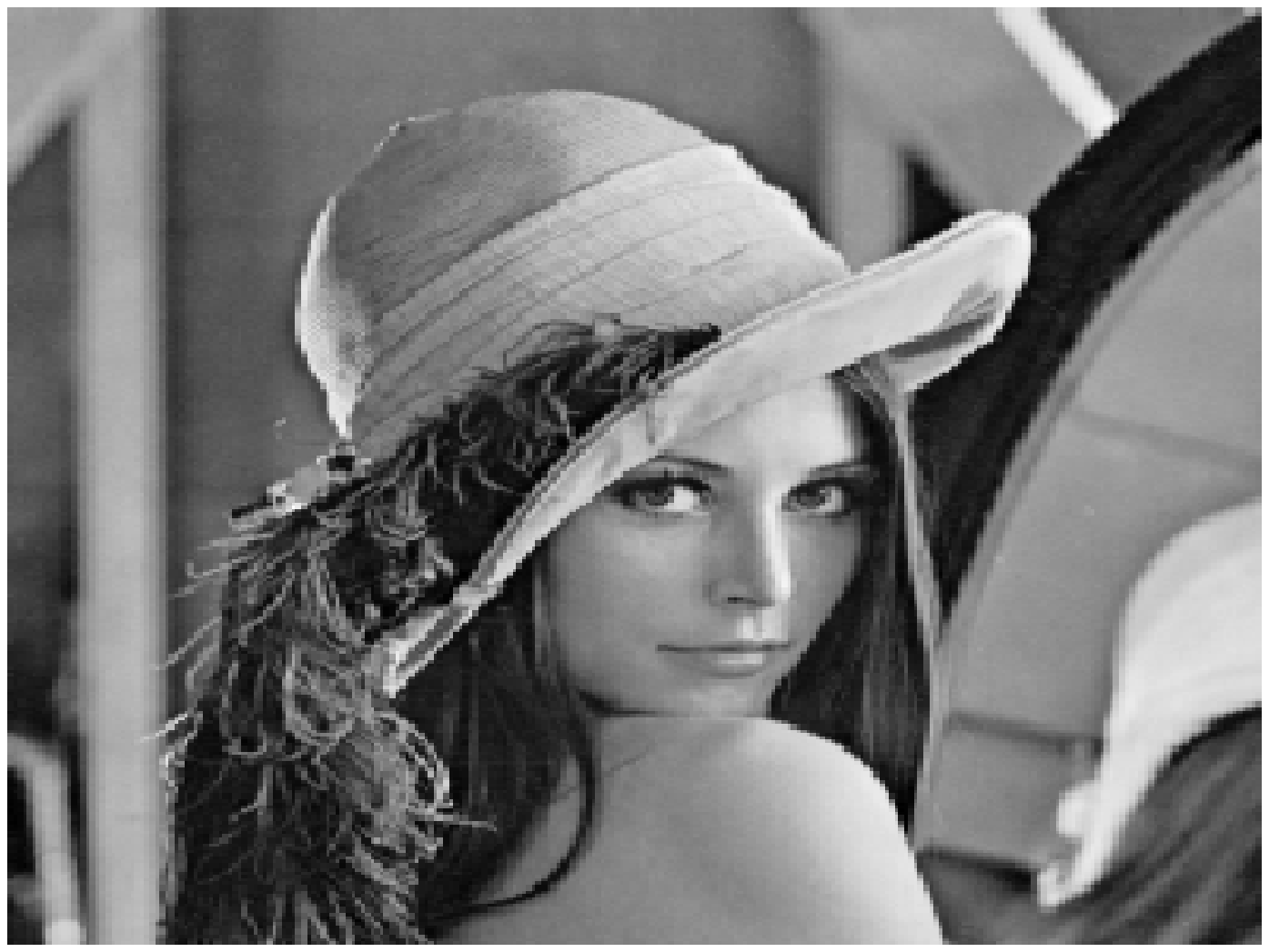} &
b)\epsfxsize=5cm \epsffile{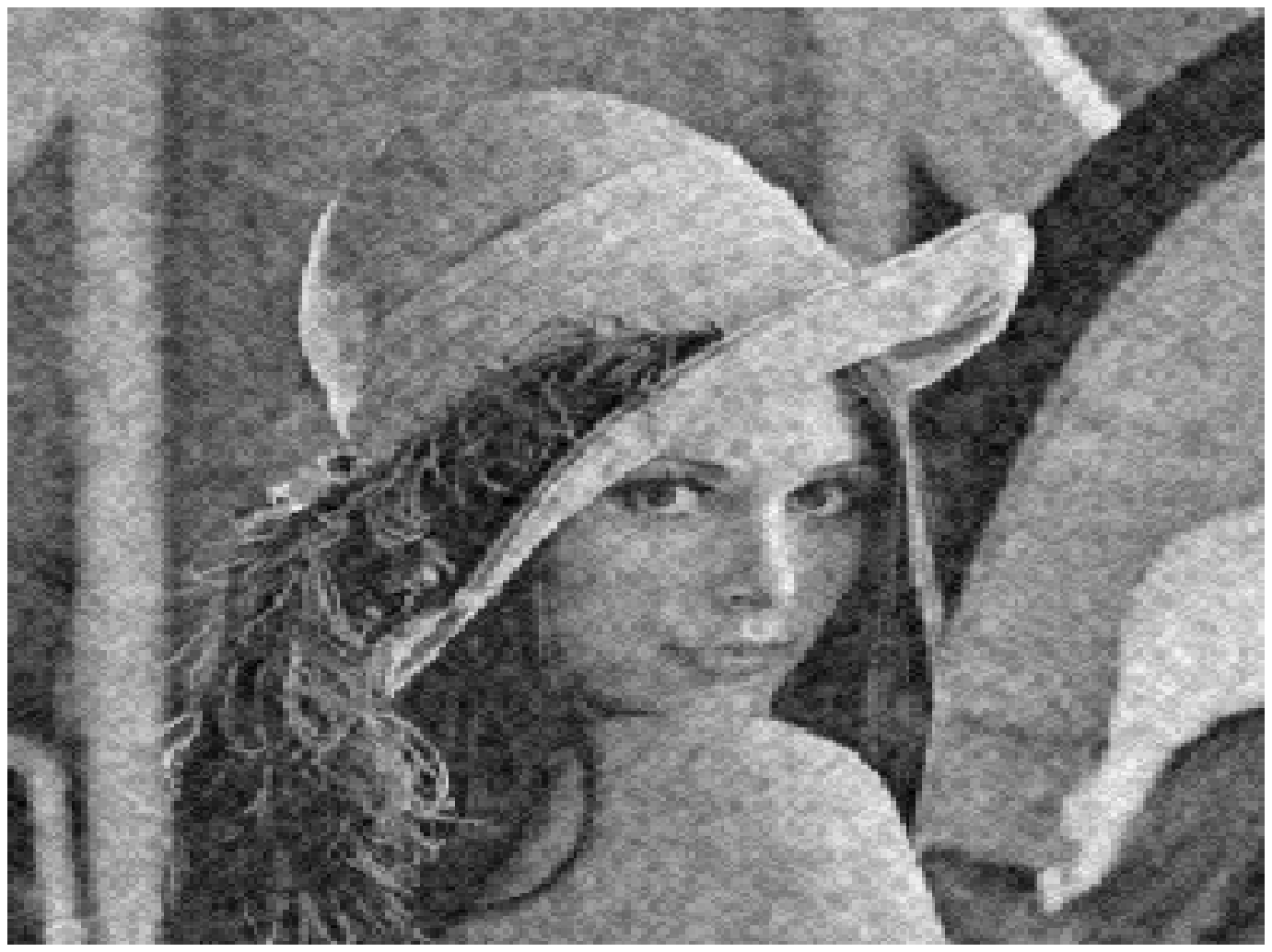}\\
c)\epsfxsize=5cm \epsffile{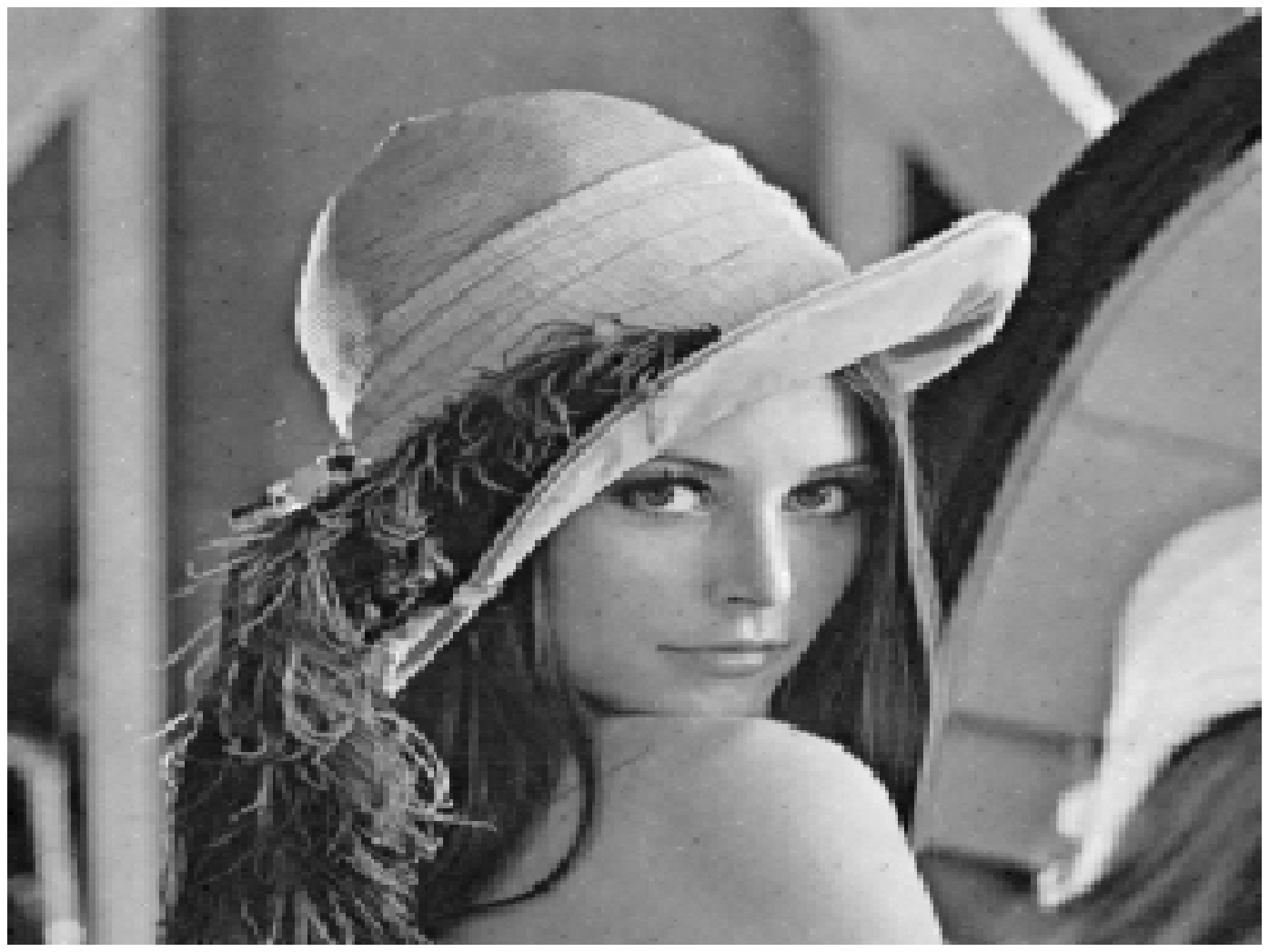} &
d)\epsfxsize=5cm \epsffile{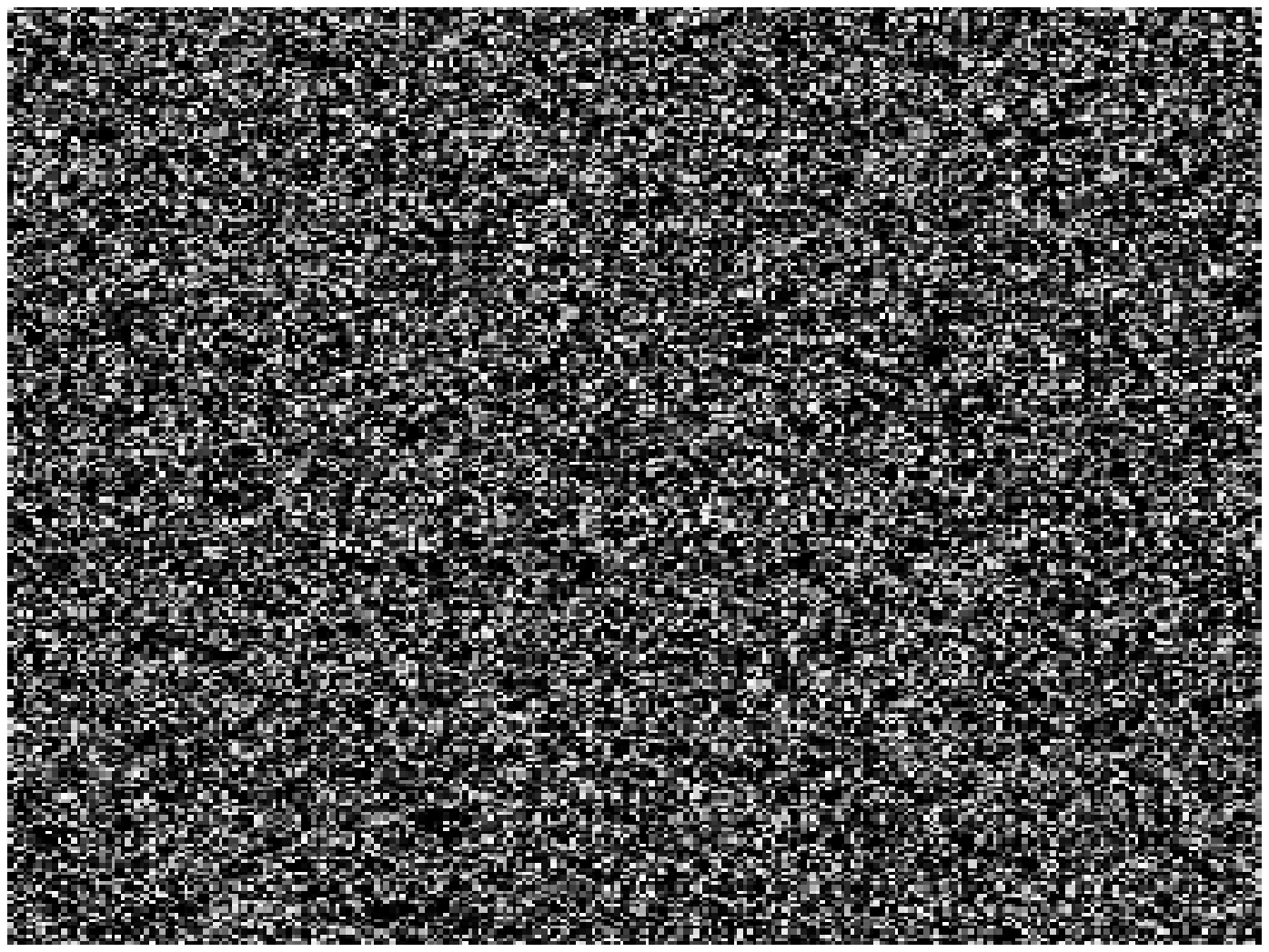}\\
\end{tabular}
\end{center}
\caption{Results of the image restoration procedure~: a) original Lena
image, b) realisation of a noisy image, c) expectation of the restored
images, d) obtained $p-$values as a gray level image (white pixels
represent values close to $1$, whereas black pixels indicate values
close to $0$).}
\label{results}
\end{figure}

%\query{... interpretation ...}
%\query{...check references ...}

% ------------------------------------------------------------------------
\bibliographystyle{plain}
\bibliography{xbib}
\vspace{+0.1cm}
\noindent
Mohamed EL MACHKOURI, Radu STOICA~\\
Laboratoire de Math\'ematiques Paul Painlev\'e\\
UMR CNRS 8524, Universit\'e Lille 1\\
U.F.R. de Math\'ematiques Pures et Appliqu\'ees\\
59655 Villeneuve d'Ascq Cedex\\
Mohamed.El-Machkouri@math.univ-lille1.fr,\\
Radu.Stoica@math.univ-lille1.fr
\end{document}